\documentclass[times,10pt]{article}

\begin{document}

\newtheorem{Remark}{Remark}[section]
\newtheorem{Lemma}{Lemma}[section]
\newtheorem{Proposition}{Proposition}[section]
\newtheorem{Theorem}{Theorem}[section]
\newtheorem{Exercise}{Exercise}[section] 
\newtheorem{Corollary}{Corollary}[section] 
\newtheorem{Definition}{Definition}[section]
\newtheorem{Example}{Example}[section]
\def\ad{\mbox{ad\,}} \def\tr{\mbox{tr\,}}\def\char{\mbox{char\,}}
\def\mod{\mbox{mod\,}}\def\Ad{\mbox{Ad\,}}
\def\diag{\mbox{diag\,}}\def\ad{\mbox{ad\,}} \def\tr{\mbox{tr\,}} 
\def\End{\mbox{End\,}}\def\GL{\mbox{GL\,}}\def\SL{\mbox{SL\,}}
\def\Out{\mbox{Out\,}}\def\Int{\mbox{Int\,}}
\def\Mat{\mbox{Mat\,}}\def\Hom{\mbox{Hom\,}}\def\Iso{\mbox{Iso\,}}
\def\Aut{\mbox{Aut\,}}\def\gl{\mbox{gl\,}}\def\sl{\mbox{sl\,}}
\def\o{\mbox{o\,}}\def\sp{\mbox{sp\,}}
\def\im{\mbox{im\,}}\def\ker{\mbox{ker\,}}\def\deg{\mbox{deg\,}}
\def\id{\mbox{id\,}}\def\mod{\mbox{mod\,}}
\def\grad{\mbox{\rm grad\,}}\def\rot{\mbox{\rm rot\,}}
\def\div{\mbox{\rm div\,}}\def\Grad{\mbox{\rm Grad\,}}
\def\det{\mbox{\rm det\,}} \def\ctg{\mbox{\rm ctg\,}}\def\tg{\mbox{\rm tg\,}}
\def\sn{\mbox{\rm sn\,}}\def\th{\mbox{\rm th\,}}\def\dn{\mbox{\rm dn\,}}
\def\pd#1,#2{\frac{\partial#1}{\partial#2}}\def\diag{\mbox{\rm diag\,}}
\def\v#1{\overline{#1}}\def\sh{\mbox{\rm sh\,}}\def\ch{\mbox{\rm ch\,}}
\def\d#1,#2{\frac{d#1}{d#2}}\def\th{\mbox{\rm th\,}}
\def\qed{\hbox{${\vcenter{\vbox{
\hrule height 0.4pt\hbox{\vrule width 0.4pt height 6pt
\kern 5pt\vrule width 0.4pt}\hrule height 0.4pt}}}$}}

\begin{titlepage} 
\title{ The Casimir Invariants of Universal Lie algebra extensions
via commutative structures}  
\author{A  B  Yanovski \thanks{\it On leave of absence from the Faculty of
Mathematics and Informatics, St.  Kliment Ohridski University, James
Boucher Blvd, 1164 - Sofia, Bulgaria}}

\maketitle 
\begin{center} 
Universidade Federal de Sergipe, Cidade Universit\'aria,\\ ''Jos\'e
Alo\'{\i}sio de Campos'', 49.100-000-S\~ao Christ\'ov\~ao, SE-Brazil
\end{center} 
\begin{abstract} 
We consider the Casimir Invariants related  to some a special kind of
Lie-algebra extensions, called universal extensions. We show that
these invariants can be studied using the equivalence between the
universal extensions and the commutative algebras and consider in
detail the so called coextension structures, arising in the
calculation of the Casimir functions.
\end{abstract}
\begin{center} 
37K30,  37K05,  17B80
\end{center} 

\end{titlepage}
\section{Introduction}
Let ${\cal G}$ be an arbitrary Lie algebra over the field ${\bf K}$.
Usually, having in mind the applications, ${\bf K}$ is considered to
be one of the classical fields ${\bf R}$ or ${\bf C}$ but there is no
immediate necessity to limit ourself. Let ${\cal G}$ be a Lie
algebra over ${\bf K}$ and let us consider the linear space ${\cal
G}^n$ with elements
\begin{equation}
{\bf x}=(x_1,x_2,\ldots,x_n),\quad x_i\in {\cal G}
\end{equation}
In \cite{JLT}, there was propounded the idea to study the 
possibilities of defining Lie brackets over ${\cal G}^n$, using the
brackets in ${\cal G}$, in such a way that the construction will be
the same for all ${\cal G}$. For this reason we call these brackets
universal ones and the corresponding Lie algebra structures universal
extensions. More specifically, let us introduce brackets obeying the
following properties:
\begin{itemize}
\item The brackets have the form:
\begin{equation}\label{eq:DEF}
([{\bf x},{\bf y}]_W)_s=\sum\limits_{i,j=1}^n W^{ij}_s[x_i,y_j];
\quad W^{ij}_s\in {\bf K}
\end{equation}
\item For $W=(W^{ij}_s)$ fixed, (\ref{eq:DEF}) is a Lie bracket for
arbitrary Lie algebra ${\cal G}$. 
\end{itemize}
For obvious reasons this algebra will be denoted by ${\cal G}^n_W$.
Let us also remark that $W=(W^{ij}_s)$ has a properties of a
$(2,1)$ tensor over ${\bf K}^n$.

The existence and classification of the above structures in the case
${\bf K}={\bf C}$ was considered in
\cite{JLT}, some additional examples of such structures are given in
\cite{ABY}. However the viewpoint in \cite{ABY} is completely
different, we tried there to establish a connection with the
so-called bundles of Lie algebras, see \cite{Ya3,TrFo}. As to the
interest in the universal extensions, it is motivated by the study of a
number of Hamiltonian structures for dynamic systems related
with Mechanics, Hydrodynamics, Magnetohydrodynamics as well as in
other areas, see
\cite{JLT} for an extensive bibliography, and the possibility to
interpret the corresponding Poisson tensors for them as Kirillov
tensors on ${\cal G}_W^*$. Recently, see \cite{JLT1}, an interesting mechanical
system having the Kirillov structure resulting from one of the most
simple extensions has been considered. (For the definitions and the
applications of the Kirillov tensors see \cite{Ki,Lich,MR}). The
structure of the universal extensions is also useful in the questions
about the stability, as there is a possibility to calculate the
Casimir invariants, see again \cite{JLT,JLT1}.

The bilinear form $[{\bf x},{\bf y}]_W$ is a Lie
bracket of the type we are looking for if the tensor $W=(W^{ij}_s)$
satisfies the conditions:
\begin{equation}\label{eq:0}
\begin{array}{c}
W^{ij}_s=W^{ji}_s\\[4pt]
\sum\limits_{k=1}^n(W_i^{sk}W_{k}^{qp}-W_i^{qk}W_{k}^{sp})=0
\end{array}
\end{equation}
In the case ${\bf K}={\bf C}$ the tensors $W_{s}^{ij}$ and the
corresponding Lie algebra extensions have been analyzed in
\cite{JLT}. It was be shown that they are divided into two classes,
called semisimple extensions and solvable extensions.  Knowing the
extensions of the second class we generate the first one and
vice-versa, so it is natural to study only the solvable extensions.
Let us briefly outline the main points of the mentioned analysis and
introduce the requisite definitions.

Instead of studying the tensors $W^{ij}_k$ it is more convenient to
study an equivalent object - the set of matrices: $W^{(i)}$, $i=1,2,\ldots, n$,
with components $(W^{(i)})^j_k=W^{ij}_k$. There is indeed a strong
reason for this, because the second equation in (\ref{eq:0})
actually means that the matrices $W^{(i)}$ commute. As a consequence
in the case ${\bf K}={\bf C}$ by a similarity transformation
$X\mapsto A^{-1}XA$ defined by a nonsingular matrix $A=(A_i^{i'})$
all the ${W}^{(s)}$ can be put simultaneously into a block-diagonal
form, each block being low-triangular with the corresponding generalized
eigenvalue on the diagonal. The linear transformation $A$ applied to
the tensor $W_k^{ij}$ of course gives
\begin{equation}\label{eq:trm1}
W^{i'j'}_{k'}=\sum\limits_{i,j,k=1}^n A^{i'}_iA^{j'}_jA^k_{k'}W_{k}^{ij}
\end{equation}
where as usual $\sum\limits_{k'}A_k^{k'}A_{k'}^s=\delta_k^s$. As a
result the set $\{W^{(i)}\}_{i=1}^n$ transforms as follows:
\begin{equation}\label{eq:trm2}
W^{(i')}=\sum\limits_{i=1}^n(A^{-1}W^{(i)}A)A^{i'}_i 
\end{equation}
and the block structure already obtained is preserved. The block
structure corresponds to a splitting of the algebra ${\cal 
G}_W^{n}$ into a direct sum. Therefore we can limit ourselves with the
irreducible case and can assume that all $W^{(s)}$ are
low triangular.  Now, the symmetry of $W_{i}^{jk}$ entails that the
generalized eigenvalues of $W^{(s)}$ for $s>1$ are zero and so for
$s>1$ the matrices $W^{(s)}$ are low-triangular with zeroes on the
diagonal, and hence are nilpotent.  All diagonal elements of
$W^{(1)}$ can be assumed to be equal to one and the same number
$a\neq 0$  or all to be equal to $0$. The first case is
called semisimple case and the second - solvable. In what follows we
shall ordinary consider solvable extensions and so for all $i$
$(W^{(i)})^j_k=0$ for $j\geq k$. In particular, $W^{(n)}$ is always
the zero matrix. In the semisimple case, as shown in \cite{JLT}, one
can make a linear transform after which $W^{(1)}={\bf 1}$, conserving
the block-triangular form of the matrices $W^{(i)}$ and the fact that
the diagonal elements of $W^{(j)}$ for $j>1$ are zero.

The lower triangular form of the matrices $W^{(i)}$
following from the above discussion in the case ${\bf K}={\bf C}$
will be called canonical form. In terms of the matrices $W^{(i)}$ the
process of passing from semisimple to solvable extension looks in the
following way. In the set $\{W^{(i)}\}_{i=1}^n$ of $n\times n$
matrices we just drop $W^{(1)}$ and form a set
$\{R^{(s)}\}_{s=1}^{n-1}$ where $R^{(s)}$ is constructed from
$W^{(s+1)}$ cutting off the first row and column.  For this reason,
in the semisimple case is useful to label the matrices by the indices
$0,1,2,\ldots n$ and in the solvable case by $1,2,\ldots, n$.  Then
$R^{(s)}$ will be obtained from $W^{(s)}$.

In \cite{Ya3} we have studied a remarkable fact - that the structure
defined by the tensor $W$ actually suggests an existence of another
structure - a structure of an associative algebra ${\cal A}_W^n$ and
there is in fact one-to-one correspondence between universal algebra
structures and associative algebra structures. In this paper we
continue to exploit this correspondence, trying to understand what
objects correspond to the Casimir functions for the Poisson-Lie
structure on $({\cal G}_W^n)^*$ on the associative algebra ${\cal A}^n_W$.
For this we shall remind the main facts, established in \cite{Ya3}.

Let ${\cal A}^n$ be a vector space with a basis
$\{e^i\}_{i=1}^n$ over the field ${\bf K}$.  Let us put:
\begin{equation}\label{eq:glt}
e^i*e^j=\sum\limits_{s=1}^n W^{ij}_s e^s
\end{equation}
where $W^{ij}_k$ is a tensor with the properties (\ref{eq:0}).
Then we have:
\begin{Theorem}\label{Theorem:1}
The $n$-dimensional universal extensions ${\cal G}^n_W$ over ${\bf
K}$ defined by the tensors $W^{ij}_k$ are in one-to-one
correspondence with the $n$-dimensional associative commutative
algebras ${\cal A}_W^n$ over the same field. Moreover, the correspondence
\begin{equation}
e^i\mapsto W^{(i)}\in \Mat(n, {\bf K}),\quad (W^{(i)})^j_k=W^{ij}_k,
\end{equation}
defines a matrix representation of the algebra ${\cal A}_W^n$ and hence
to any linear transformation of the elements $e^i$
\begin{equation}
e^{i'}=\sum\limits_{i=1}^n A^{i'}_i e^i
\end{equation}
correspond the transformations (\ref{eq:trm1}) and (\ref{eq:trm2})
of the tensor $W^{ij}_k$ and of the matrices $W^{(i)}$.
\end{Theorem}
From the above result it is patent that the matrix $W^{(i)}$ is in
fact the matrix of the action of the element $e^i$ on ${\cal A}_W^n$.
This permits to understand easily the block structure and the
canonical structure of the matrices $W^{(i)}$:
\begin{Proposition}
The splitting of ${\cal A}_W^n$ into a sum of ideals:
\begin{equation}
{\cal A}_W^n={\cal J}^{s_1}_{W_1}\oplus {\cal J}^{s_2}_{W_2}\oplus
\ldots\oplus{\cal J}^{s_q}_{W_q} 
\end{equation}
with dimensions $s_i=\dim{\cal J}^{s_i}_{W_i}$ corresponds to the
splitting of the algebra ${\cal G}_W$ into a direct sum:
\begin{equation}
{\cal G}_W^n={\cal G}^{s_1}_{W_1}\oplus{\cal
G}^{s_2}_{W_2}\oplus\ldots \oplus{\cal G}^{s_q}_{W_q}
\end{equation}
The irreducible ${\cal G}^{n}_{W}$ corresponds to ${\cal A}_W^n$ that
cannot be split into a sum of proper ideals, that is to simple
${\cal A}_W^n$. 
\end{Proposition}
\begin{Proposition}
The canonical structure of $W^{(i)}$ is equivalent to the requirement
that for the solvable case in the basis $\{e^i\}_{i=1}^n$  
\begin{equation}
\begin{array}{lr}
e^i*e^j=
\sum\limits_{s>max(i,j)}W^{ij}_s e^s&1\leq i,j\leq n
\end{array}
\end{equation}
and in the semisimple case we have one more independent element $e^0$
such that 
\begin{equation}
\begin{array}{lr}
e^0*e^j=e^j,&0\leq j\leq n\\
\end{array}
\end{equation}
\end{Proposition}
Thus the solvable case corresponds to the situation when
the operators $T_a$, $a\in {\cal A}^n_W$
\begin{equation}
T_a(b)\equiv a*b
\end{equation}
are nilpotent and the semisimple case corresponds to algebras with unity.

\section{Casimir invariants}

As already mentioned, one of the motivations in considering the
universal extensions is their applications to the  theory of Poisson
structures. It is well known that on the dual ${\cal H}^*$ of a Lie
algebra ${\cal H}$ over ${\bf K}$, provide one can identify the
spaces $[{\cal H}^*]^*$ and ${\cal H}^*$, one can define in a
canonical way Poisson bracket. Thus the space of the smooth functions
${\cal M}\equiv C^{\infty}({\cal H}^*)$ over ${\cal H}^*$ becomes a
Lie algebra. The bracket is known as Poisson-Lie bracket and the
corresponding tensor as Kirillov tensor, \cite{Ki,MR,TrFo}. The
construction is the following.  If $f,g\in {\cal M}$ and $\mu\in
{\cal H}^*$ then the differentials $df|_{\mu}, dg|_{\mu}$ can be
considered as elements from ${\cal H}$ and we can set:
\begin{equation}
\{f,g\}(\mu)=\langle [df|_{\mu}, dg|_{\mu}],\mu\rangle_{\cal H},
\end{equation}
where by $\langle~,~\rangle_{\cal H}$ is denote the canonical pairing
between ${\cal H}$ and ${\cal H}^*$. One can also write:
\begin{equation}
\{f,g\}(\mu)=\langle dg|_{\mu}, \ad^*_{df|_{\mu}}(\mu)\rangle_{\cal H},
\end{equation}
where $x\mapsto \ad_x\in \End({\cal H})$ is the adjoint
representation of ${\cal H}$ ($\ad_x(y)=[x,y]$, $x,y\in {\cal H}$) and 
\begin{equation}
\begin{array}{c}
x\mapsto-\ad^*_x\in \End({\cal H}^*)\\[4pt]
\langle y,\ad^*_{x}(\mu)\rangle_{\cal H}=\langle
\ad_x(y),\mu\rangle_{\cal H}\\[4pt]
x,y\in {\cal H},\quad \mu\in {\cal H}^*
\end{array}
\end{equation}
is the coadjoint representation. By definition the Casimir functions
are those that belong to the center of the algebra ${\cal M}$. In the
end of these general remarks we remind, that when the algebra ${\cal
H}$ splits into a direct sum, the dual also splits:
\begin{equation}
{\cal H}^*={\cal H}^*_1\oplus{\cal H}^*_2\oplus\ldots\oplus {\cal H}^*_k
\end{equation}
Then if $\pi_i$ is the projection onto ${\cal H}^*_i$ and $C_i$ is a
Casimir function on ${\cal H}^*_i$ the function $C_i\circ \pi_i$ is a
Casimir function for ${\cal H}^*$. Conversely, if $C$ is is a
Casimir function for ${\cal H}^*$ we obtain a Casimir function for
${\cal H}_i^*$ restricting it to ${\cal H}_i^*$. Thus usually only
the irreducible case is studied.

We shall apply now the above construction
for the algebra ${\cal G}_W^n$. The question about the Casimir functions 
for it is considered in detail in \cite{JLT}. However, according to our
ideology, here we shall try to describe the Casimir functions
{\it only in terms of the associative algebra ${\cal A}^n_W$}.

Let us denote by $\langle~,~\rangle_{\cal G}$ the pairing between
${\cal G}$ and ${\cal G}^*$. Next, denote by ${\xi},{\eta},\ldots$
the elements from $\left({\cal G}_W^n\right)^*=({\cal G}^*)^n$:
\begin{equation}
{\xi}=(\xi^1,\xi^2,\ldots,\xi^n), \quad \xi_j\in {\cal G}^*
\end{equation}
and the pairing between ${\cal G}^n_W$ and $({\cal
G}^n)^*$ by $\langle~,~\rangle_W$. Now, suppose that $f\in {\cal
M}_W=C^{\infty}(({\cal G}_W^*)^n)$ and let us denote by $\partial_i
f|_{\xi}=(\partial_i f)_{\xi}$ the partial derivative with
respect to $\xi_i$, that is:
\begin{equation}
(\partial_i f)_{\xi}(\eta_i)=df|_{\xi}(0,0,...,\eta_i,...,0,0), \quad
\eta_i\in {\cal G}^*
\end{equation}
Naturally, we consider $(\partial_i f)_{\xi}$ as elements from
${\cal G}$. Thus finally, we can write the Poisson-Lie bracket over
$\left({\cal G}^n_W\right)^*$. For $f,g\in {\cal
M}_W$
\begin{equation}
\{f,g\}_W=\sum\limits_{i,j,k=1}^nW^{ij}_k\langle [\partial_i f,
\partial_j g],\xi^k\rangle_{\cal G}
\end{equation}
The Casimir function $C({\xi})\in {\cal M}_W$ therefore must
satisfy the equation:
\begin{equation}\label{eq:Ca}
\sum\limits_{i,j,k=1}^nW^{ij}_k\langle [\partial_i C,
\partial_j f],\xi^k\rangle_{\cal G}=0, \quad f\in {\cal M}_W
\end{equation}
or simply
\begin{equation}\label{eq:CCON}
\sum\limits_{i,j=1}^nW^{ij}_k\ad^*_{\partial_i C}(\xi^k)=0
\end{equation}
The most interesting cases arise when ${\cal G}$ and ${\cal G}^*$
can be canonically identified through a nondegenerate invariant
bilinear form. If ${\cal G}$ is finite dimensional, the existence of
such a form in the case ${\bf K}={\bf R}$ or ${\bf K}={\bf C}$ means
that the algebra is semisimple and usually as such a form is taken
the Killing form
\begin{equation}
K(x,y)=\tr(ad_x\circ\ad_y),\quad x,y\in {\cal G}
\end{equation}
When the algebra is simple the invariant nondegenerate form is defined
up to a constant multiplier, but in the semisimple case it can be a
linear combinations of the forms for the simple components, see \cite{GoGr}.
In any case, we shall assume that such a form exists. Instead of
$K(x,y)$ we shall prefer to work with the bijective linear map
\begin{equation}
\begin{array}{c}
\hat{K}:{\cal G}\mapsto {\cal G}^*\\[4pt]
\langle x,\hat{K}y\rangle_{\cal G}=K(x,y),\quad x,y\in {\cal G}
\end{array}
\end{equation}
The invariance of $K$ means that
\begin{equation}
\hat{K}\circ \ad_x=-\ad_x^*\circ \hat{K}
\end{equation}
and then the adjoin and coadjoint representations are equivalent.
Below we shall consider the case when ${\cal G}$ and ${\cal G}^*$ can
be identified through $\hat{K}$, that is we limit ourself to the
semisimple ${\cal G}$.

Returning to the question of the Casimir functions it is natural
to consider polynomial $C$ and therefore we must start
with homogeneous functions. The most simple case is clearly
a linear function:
\begin{equation}
C^{(1)}({\xi})=\langle {y}, {\xi}\rangle_W=\sum\limits_{i=1}^n
\langle y_i,\xi^i\rangle_{\cal G}, \quad {y}\in {\cal G}
\end{equation}
The equation (\ref{eq:Ca}) then can be written as:
\begin{equation}
\sum\limits_{i,j,k=1}^n\langle W_k^{ij}y_i,\ad^*_{\partial_j
f}\xi^k\rangle_{\cal G}=0
\end{equation}
Setting $x^k=\hat{K}^{-1}\xi^k$ we put the above into the form:
\begin{equation}
\sum\limits_{i,j,k=1}^n\langle \hat{K}(W_k^{ij}y_i),[{\partial_j
f}, x^k]\rangle_{\cal G}=0
\end{equation}
and if we want the above to be satisfied for arbitrary semisimple algebra
${\cal G}$ and arbitrary $f$ we obtain
\begin{equation}
\sum\limits_{i=1}^n W_k^{ij}y_i=0,\quad 1\leq k,j\leq n
\end{equation}
This is obviously true for $y_i=p_iy$, where $y\in {\cal G}$ and
$p_i\in {\bf K}$, such that the vector $({p})_i=p_i$ from ${\bf
K}^n$ is eigenvector with zero eigenvalue for all the matrices
$W^{(j)}$. For example, if $W^{(j)}$ are in canonical form, then
$p_i=\delta_i^n$ is an eigenvector with zero eigenvalue. 
The vector ${p}$ defines an element from ${\cal A}^n_W$,
$p=\sum_{i=1}^n p_ie^i$.  Then we have:
\begin{Proposition}
To arbitrary linear function on $({\cal A}^n_W)^*$ of the type:
\begin{equation}
\begin{array}{c}
\alpha\mapsto \langle p,\alpha\rangle_{\cal A}\\[4pt]
p\in {\cal A}^n_W,~ T_b(p)=b*p=0,~ b\in {\cal A}_W^n
\end{array}
\end{equation}
corresponds Casimir function of the type $C^{(1)}$.
\end{Proposition}

As we have seen in the so called solvable case such vectors exist and in
the so called semisimple case they do not exist.

The second simple class is of course the class of quadratic Casimirs.
We shall suppose that  they are of the type:
\begin{equation}
C^{(2)}={1\over 2}\sum\limits_{i,j=1}^n C_{ij} \langle
\hat{K}^{-1}\xi^i,\xi^j\rangle_{\cal G}
\end{equation}
where $C_{ij}=C_{ji}$ is a symmetric matrix. From (\ref{eq:CCON})
we easily get the condition defining the Casimirs:
\begin{equation}\label{eq:Ca2}
\sum\limits_{i=1}^n(W_k^{ji}C_{ir}-W_{r}^{ji}C_{ik})=0,\quad
0\leq j,k,r\leq n
\end{equation}
If we now set 
\begin{equation}
C(x,y)=\sum\limits_{i,j=1}^nC(e_i,e_j)x^iy^j\equiv
\sum\limits_{i,j=1}^nC_{ij}x^iy^j
\end{equation}
for $x=\sum_{i=1}^nx^ie_i,~y=\sum_{i=1}^nx^ie_i\in ({\cal A}^n_W)^*$
then we obtain the following:
\begin{Proposition}
To every symmetric
bilinear form $C(x,y)$ on $({\cal A}^n_W)^*$, such that
\begin{equation}
C(T_a^*x,y)-C(x,T^*_ay)=0,\quad a\in {\cal A}^n_W, \quad x,y\in
({\cal A}^n_W)^* 
\end{equation}
corresponds a Casimir function of the type $C^{(2)}$.
\end{Proposition}
Following our idea we shall call the above functions Casimir
functions of ${\cal A}_W^n$, though of course this means that they
are Casimir functions for ${\cal G}^{n}_{W}$.
\begin{Example}
A fact first noted in \cite{JLT} is that a special solution of
(\ref{eq:Ca2}) exists if the vector ${n}\in {\bf K}^n$  is an
eigenvector for all the matrices $W^{(j)}$. This solution is 
$C_{ij}=n_in_j$. In covariant notations, to ${n}$ corresponds an
element $n=\sum_{i=1}^n n_ie^i\in {\cal A}^n_W$ and to it a bilinear
form over (${\cal A}^n_W)^*$
\begin{equation}
C(x,y)=\langle n\otimes n, x\otimes y\rangle_{\cal A} = \langle n,
x\rangle_{\cal A}\langle n,y\rangle_{\cal A}
\end{equation}
\end{Example}
Oddly enough, the equation (\ref{eq:Ca2}) characterizes the Casimirs in
general, in the case when the algebra ${\cal G}$ is the infinite dimensional
algebra of the real-valued functions $f$ with support laying in an
open connected region $\Omega$ of the real plane ${\bf R}^2$ with the
bracket:
\begin{equation}
[f,g](x,y)={\partial f\over \partial x}{\partial g\over \partial y}-
{\partial f\over \partial y}{\partial g\over \partial y},\quad
(x,y)\in \Omega
\end{equation}
see \cite{JLT}. The corresponding Hamiltonian structures have
important applications in the hydrodynamics, so the study of the
equations (\ref{eq:Ca2}) is very interesting. The algebras ${\cal G}$
and ${\cal G}^*$ are identified through the bilinear form:
\begin{equation}
(f,g)=\int\limits_{\Omega}fg~dxdy
\end{equation}
and the Casimirs are functionals of the form:
\begin{equation}\label{eq:VID}
\hat{C}[{\bf \xi}]=\int\limits_{\Omega}C(\xi^1,\xi^2,\ldots,\xi^n)dx dy
\end{equation}
$\xi^i=\xi^i(x,y)\in {\cal G}^*$. The function
$C=C(\xi^1,\xi^2,\ldots,\xi^n)$ is a "usual" differentiable function.
Then the calculations show that the equations for the Casimir
function of the type (\ref{eq:VID}) are exactly (\ref{eq:Ca2}), where
$C_{ij}={\partial^2 C\over \partial \xi^i\partial\xi^j}$.

\subsection{Reduction of the Casimir functions}

We shall try now to describe the reduction of the Casimir functions,
outlined in \cite{JLT} on the base of equations (\ref{eq:Ca2})
within the new framework, that of the associative algebra ${\cal
A}^n_W$.  As we shall see, we do not need the existence of the
canonical base to perform the main construction. To this end we shall
need only the requirements that in the algebra ${\cal A}^{n+1}_W$
there is a basis $\{e^i\}_{i=0}^n$ such that $e^0$ is unity, that is
$e^0*e^i=e^i$, $0\leq i\leq n$ and as we shall say there exists a pseudo-zero
element, that is $e^n$, such that $e^n*e^i=\delta^{i}_{0}e^n$, $0\leq
i\leq n$. In other words, we assume:
\begin{equation}\label{eq:AS}
W^{0i}_j=\delta_{j}^i,\quad W^{ni}_j=\delta_0^i\delta_{j}^{n}
\end{equation}
It is readily seen that from (\ref{eq:0}) follows that $W^{ij}_0=0$,
$1\leq i,j\leq n-1$, a fact that we shall use later.

In the semisimple case provide we are working in a canonical base
$\{e^i\}_0^n$ the vectors $e^0$ and $e^n$ possess the above
properties, in the solvable case in a canonical base $\{e_i\}_1^n$
one must just add the ''solvable part'', that is a unity $e_0$ in
order to come to the situation we need. Our conditions 
are a little weaker tan those in \cite{JLT}, we just need the
properties of the first and last elements of the canonical base. 

We can express also (\ref{eq:AS}) in terms of the adjoint action:
\begin{equation}
T_{e^0}^*e_i=e_i,\quad T_{e^n}^*e_i=\delta^n_ie_0
\end{equation}
In particular, $T_{e^n}^*e_n=e_0$.
The invariance conditions on the components $C_{ii}$ are trivially
satisfied, that is why we shall consider only the case $C_{ij}$,
$i\neq j$. We note also, that as $T^*_{e^0}=\id$ the invariance with
respect to this transform is trivial. Therefore it is enough to
consider invariance relations with respect to $T^*_{e^j}$, $j\neq 0$.
Let us start with $T^*_{e^n}$. On the first place we can prove that
$C_{0i}=C(e_0,e_i)=0$, $1\leq i\leq n-1$. Indeed, $$
C(e_0,e_i)=C(T^*_{e^n}e_n,e_i)=C(e_n,T^*_{e^n}e_i)=\delta_i^nC(e_0,e_n)=0
$$ 
We have also $C_{00}=0$, for 
$$
0=C(T^*_{e^n}e_0,e_n)=C(e_0,T^*_{e^n}e_n)=C(e_0,e_0) $$ 
So $C_{0j}=0$; $0\leq j\leq n-1$ and we can see that the invariance
equations corresponding to $T^*_{e^n}$ are satisfied automatically.
Indeed, for $1\leq i\leq n-1$ we have: 
$$
\begin{array}{c}
C(T^*_{e^n}e_i,e_n)=\delta_i^n C(e_0,e_n)=0\\[4pt]
C(e_i,T^*_{e^n}e_n)=C(e_i,e_0)=0
\end{array}
$$
and hence always $C(T^*_{e^n}e_i,e_n)=C(e_i,T^*_{e^n}e_n)$. Also, it is
readily seen that 
$$
\begin{array}{c}
C(T^*_{e^n}e_n,e_0)=C(e_n,T^*_{e^n}e_0)=0\\[4pt]
C(T^*_{e^n}e_i,e_j)=C(e_i,T^*_{e^n}e_j)=0\\[4pt]
1\leq i,j\leq n-1
\end{array}
$$

Let us pass now to the invariance conditions corresponding to
$T^*_{e^j}$, $j\neq 0,n$. For $1\leq s,i,j\leq n-1$ we have 
$$C(T^*_{e^s}e_i,e_j)=\sum\limits_{k=1}^{n-1}W_{i}^{sk}C(e_k,e_j)+
W_i^{s0}C(e_0,e_j)+W_i^{ns}C(e_n,e_j) $$ The last two terms being
equal to zero, we obtain
\begin{equation}
C(T^*_{e^s}e_i,e_j)=\sum\limits_{k=1}^{n-1}W_{i}^{sk}C(e_k,e_j)
\end{equation}
so  the invariance conditions for $C_{ij}$, $1\leq i,j\leq n-1$ with
respect to $T^*_{e^s}$, $1\leq s\leq n-1$ are written only in terms
of values having indices running from $1$ to $n-1$:
\begin{equation}\label{eq:Caz1}
\begin{array}{c}
\sum\limits_{i=1}^{n-1}(W_k^{ji}C_{ir}-W_{r}^{ji}C_{ik})=0\\[4pt]
1\leq j,k,r\leq n-1
\end{array}
\end{equation}
As to the coefficients $C(e_i,e_n)$, it is easy to see that the
invariance condition:
$C(T^*_{e^s}e_i,e_n)=C(e_i,T^*_{e^s}e_n)$, $1\leq
i,s\leq n-1$ is equivalent to:
\begin{equation}\label{eq:Caz2}
\begin{array}{c}
\sum\limits_{m=1}^{n-1}W_i^{sm}C(e_m,e_n)-\sum\limits_{m=1}^{n-1}W_n^{sm}C(e_m,e_i)
+\delta_{i}^sC(e_0,e_n)=0\\[4pt]
1\leq i,s\leq n-1
\end{array}
\end{equation}
Let the the solvable extension corresponding to ${\cal A}^{n+1}_{W}$ 
be ${\cal A}^n_{\hat{W}}$ (We identify the extensions with the 
corresponding associative algebras).  ${\cal A}^n_{\hat{W}}$ is
obtained putting formally $e^0=0$.  In ${\cal A}^n_{\hat{W}}$ the space
${\bf K}e^n$ is an ideal, as multiplication by $e_n$ always gives
zero.  Factoring over it we obtain an algebra ${\cal
A}^{n-1}_{\bar{W}}$, $\bar{W}^{ij}_k=W^{ij}_k$, $1\leq i,j,k\leq
n-1$. It is clear, that $C$ can be restricted to ${\cal
A}^{n-1}_{\bar{W}}$. We shall denote this restriction with the letter
$\bar{C}$. More generally, if we have some object related to ${\cal
A}^{n+1}_W$ and we denote it say by $A$, then the analogous object
related with ${\cal A}^n_{\hat{W}}$ we shall denote by $\hat{A}$ and the
corresponding object related with ${\cal A}^{n-1}_{\bar{W}}$ by
$\bar{A}$. Going back to the question of the Casimir function
$\bar{C}$ we get:
\begin{Proposition}[see \cite{JLT}]
The invariance of $C$ with respect of the adjoint action of ${\cal
A}^{n+1}_W$ is equivalent to the invariance of $\bar{C}$ with respect to the 
adjoint action of ${\cal A}^{n-1}_{\bar{W}}$ and the condition
(\ref{eq:Caz2}) over the coefficients. 
\end{Proposition}
\begin{Corollary}
The function $\bar{C}$ is a Casimir function for the algebra ${\cal
G}^{n-1}_{\bar{W}}$. 
\end{Corollary}
If ${\cal A}^{n-1}_{\bar{W}}$ is an arbitrary associative algebra,
adding a pseudo-zero and unity we can obtain algebra ${\cal
A}^{n+1}_{W}$ of the type we started with. Thus we have:
\begin{Corollary}
If $C_{ij}$, $1\leq i,j\leq n-1$ is a Casimir function for ${\cal
A}^{n-1}_{\bar{W}}$ then for arbitrary solution $X_i$, $0\leq i\leq
n-1$ of the system:
\begin{equation}
\begin{array}{c}
\delta_{i}^sX_0+\sum\limits_{m=1}^{n-1}W_i^{sm}X_m=
\sum\limits_{m=1}^{n-1}W_n^{sm}C_{mi}\\
1\leq i,s\leq n-1
\end{array}
\end{equation}
and arbitrary number $C_{nn}$ we can obtain a Casimir
function for ${\cal A}^{n+1}_{W}$ setting $C(e_n,e_n)=C_{nn}$,
$C(e_0,e_i)=0$, $1\leq i\leq n-1$ and $C(e_n,e_j)=X_j$, $0\leq j\leq n-1$.
\end{Corollary}

\subsection{Calculation of the Casimirs and Coextensions}

Looking upon the system (\ref{eq:Caz2}) it can be seen that we can
consider as given the quantities $C_{in}$ and try to resolve for
$C_{ij}$, $1\leq i,j\leq n-1$, that is, we can try to construct
Casimirs for the higher dimensions. On this consideration is based 
the ingenious approach used in \cite{JLT} in order to construct the
Casimirs. Below we shall give coordinate-free description of this approach
and to the so-call coextension which plays crucial role in it. According to
\cite{JLT} there are two distinct cases, which must be considered
separately -  when the matrix $(W_n^{ij})_{1\leq i,j\leq n-1}$ is
nondegenerate and when it is degenerate. We shall concentrate upon
the first possibility, the results about the degenerate case will be
presented elsewhere. In coordinate free notation, the fact that $(W_n^{ij})_{1\leq i,j\leq n-1}$ is nondegenerate means that the map:
\begin{equation}\label{eq:PSI}
\begin{array}{c}
\bar{\Psi}:{\cal A}_{\bar{W}}^{n-1}\mapsto ({\cal A}_{\bar{W}}^{n-1})^*\\[4pt]
\bar{\Psi}(x)=T^*_{x}e_n
\end{array}
\end{equation}
is an isomorphism. (Here and below we use the convention accepted
earlier for the quantities related with the tensors $W,\hat{W}$ and $\bar{W}$.)

We have
\begin{Proposition}
If the map (\ref{eq:PSI}) is isomorphism, then the invariance conditions:
$$C(T^*_{e^s}e_i,e_j)=C(e_i,T^*_{e^s}e_j)$$ follow from the
conditions $$C(T^*_{e^s}e_n,e_j)=C(e_n,T^*_{e^s}e_j)$$
($1\leq i,j,s\leq n-1$).
\end{Proposition}
Indeed, let us fix $i$. Then according to our assumption there exists
$f_i\in {\cal A}^{n-1}_{\bar{W}}$ such that $T^*_{f_i}e_n=e_i$.
We can write:
$$
\begin{array}{c}
C(T^*_{e^s}e_i,e_j)=C(T^*_{e^s}T^*_{f_i}e_n,e_j)=
C(T^*_{e^s*f_i}e_n,e_j)=C(e_n,T^*_{e^s*f_i}e_j)=\\[8pt]
C(e_n,T^*_{f_i*e^s}e_j)=C(T^*_{f_i}e_n,T^*_{e^s}e_j)=
C(e_i,T^*_{e^s}e_j)
\end{array}
$$
The coordinate form of this proposition is that provide
$(W_{n}^{ij})_{1\leq i,j\leq n-1}$is nondegenerate, the equations
(\ref{eq:Caz1}) follow from equations (\ref{eq:Caz2}), see \cite{JLT}.

The map $\bar{\Psi}$ can be used to transfer the algebra structure from
${\cal A}_{\bar{W}}^{n-1}$ onto the dual space $({\cal A}_{\bar{W}}^{n-1})^*$.

\begin{Proposition}\label{Prop:GLAI}
If the map $\bar{\Psi}$ defined in (\ref{eq:PSI})  is isomorphism, then
there exists a unique algebra structure over $({\cal
A}^{n-1}_{\bar{W}})^*$ such, that $\bar{\Psi}$ is an isomorphism between
algebras ${\cal A}_{\bar{W}}^{n-1}$ and $({\cal
A}_{\bar{W}}^{n-1})^*$.
\end{Proposition}
{\it Proof.} Let us denote by $\bar{*}$ the multiplication in ${\cal
A}^{n-1}_{\bar{W}}$ and by $\bar{T}_a$ the corresponding action
$\bar{T}_a b=a\bar{*}b$. Then if $\xi,\eta \in ({\cal
A}_{\bar{W}}^{n-1})^*$ we can define 
\begin{equation}
\xi\bar{*}\eta=\bar{\Psi}(\bar{\Psi}^{-1}(\xi)\bar{*}\bar{\Psi}^{-1}(\eta))
\end{equation}
Thus the bilinearity, the commutativity and the associativity as well
as the uniqueness are trivial to prove. Q.E.D.

It is not difficult to prove also that
\begin{equation}\label{eq:stno1}
\xi\bar{*}\eta=\bar{T}^*_{\bar{\Psi}^{-1}(\eta)}\xi
\end{equation}
and
\begin{equation}\label{eq:stno}
T_{x}^*\xi=\bar{T}_{x}^*\xi+\langle x,\xi\rangle e_0,\quad 
x\in {\cal A}_{\bar{W}}^{n-1}, \xi\in ({\cal A}_{\bar{W}}^{n-1})^*
\end{equation}
where $\bar{T}^*$ is the coadjoint action corresponding to $\bar{*}$.

The above proposition gives coordinate-free description in terms of
the commutative structures of the "coextension", introduced in
(\cite{JLT}). It is defined there as a tensor $\bar{A}_{ij}^k$
\begin{equation}
\bar{A}^i_{jk}=\sum\limits_{s=1}^{n-1}\bar{g}_{js}\bar{W}^{si}_k
\end{equation}
where $\bar{g}_{ij}$ are the components of the matrix, inverse to
$(W^{ij}_n)_{1\leq i,j\leq n-1}$.

As we can see it corresponds to the algebra $({\cal
A}_{\bar{W}}^{n-1})^*$ in the sense that
\begin{equation}
e_i\bar{*} e_j=\sum\limits_{k=1}^{n-1}\bar{A}_{ij}^ke_k
\end{equation}

Now it is easy to calculate $C_{ij}$, $1\leq i,j\leq n-1$. For $e_i$
there exists unique $f_i$ such that $\bar{\Psi}(f_i)=e_i$, actually:
\begin{equation}
f_i=\sum\limits_{k=1}^{n-1}\bar{g}_{ik}e^k
\end{equation}
Then using (\ref{eq:stno1}),(\ref{eq:stno}) we get
$$C_{ij}=C(e_i,e_j)=C(T^*_{f_i}e_n,e_j)=C(e_n,T^*_{f_i}e_j)$$
In other words:
\begin{equation}\label{eq:CC1}
C_{ij}=C(e_n,T^*_{{\Psi}^{-1}(e_i)}e_j)=C(e_n,e_i\bar{*}
e_j)+C(e_n,e_0)\bar{g}_{ij} 
\end{equation}
or more explicitly:
\begin{equation}\label{eq:CC2}
C_{ij}=\bar{g}_{ij}C_{0n}+\sum\limits_{k=1}^{n-1}\bar{A}_{ij}^k C_{nk}
\end{equation}

However, actually all this is a part of the whole algebraic picture,
a fact that seems overlooked until now. The point is that the map
$\bar{\Psi}$ can be extended to act onto the space ${\cal A}_W^{n+1}$
by the same formula:
\begin{equation}
\Psi(x)=T^*_xe_n
\end{equation}
(we denote the extension by the letter $\Psi$). Then
\begin{equation}
\Psi(e^n)=T^*_{e^n}e_n=e_0,\quad \Psi(e^0)=T^*_{e^0}e_n=e_n
\end{equation}
the space $V$ spanned by $\{e^0,e^n\}$ is mapped by
$\Psi$ into $V^*$ and onto $D={\cal A}^{n-1}_{\bar{W}}$ we have
$\Psi|_{D}=\bar{\Psi}$. The map $\Psi$ is nondegenerate exactly when
$\bar{\Psi}$ is nondegenerate. Now $\Psi$ can be used to define an
algebra structure over $({\cal A}_W^{n+1})^*$:
\begin{equation}
\xi*\eta=\Psi(\Psi^{-1}(\xi)*\Psi^{-1}(\eta)), \quad \xi,\eta
\in ({\cal A}_W^{n+1})^*
\end{equation}
If we write down the elements $\xi,\eta$ in the basis $\{e_i\}_{0}^n$
\begin{equation}
\begin{array}{c}
\xi=\xi^0e_0+\sum\limits_{i=1}^{n-1}\xi^i e_i+\xi^ne_n\\[4pt]
\eta=\eta^0e_0+\sum\limits_{i=1}^{n-1}\eta^ie_i+\eta^ne_n
\end{array}
\end{equation}
it is not difficult to calculate that
\begin{equation}\label{eq:gens}
\begin{array}{c}
\xi*\eta=(\sum\limits_{i,j=1}^{n-1}\bar{g}_{ij}\xi^i\eta^j +\xi^n\eta^0+
\eta^n\xi^0)e_0+\\[4pt]
\sum\limits_{k=1}^{n-1}(
\sum\limits_{i,j=1}^{n-1}\xi^i\eta^j\bar{A}_{ij}^k+\xi^n\eta^k+\eta^n\xi^k)e_k+
(\xi^n\eta^n)e_n
\end{array}
\end{equation}
In the above calculation we have used here the fact that 
$W_0^{ij}=0$, $1\leq i,j\leq n-1$, also $\bar{g}_{ij}$ and
$\bar{A}_{ij}^k$ have the same meaning as before.

It natural to call this structure a coextension, the one we had
formerly, as we shall see, is a consequence from this one. Is
interesting to note that for the coextension $e_0$ now is
pseudo-zero, and $e_n$ is the unity, as
\begin{equation}
e_0*\eta=\eta^ne_0,\quad e_n*\eta=\eta
\end{equation}
The solvable part, as usual, is obtain putting formally $e_n=0$ and
$\xi^n=\eta^n=0$: 
\begin{equation}
\begin{array}{c}
\xi\hat{*}\eta=(\sum\limits_{i,j=1}^{n-1}\bar{g}_{ij}\xi^i\eta^j)e_0+
\sum\limits_{k=1}^{n-1}(\sum\limits_{i,j=1}^{n-1}\xi^i\eta^j
\bar{A}_{ij}^k)e_k\\[8pt]
\xi,\eta\in ({\cal A}^{n}_{\hat{W}})^*
\end{array}
\end{equation}
or
\begin{equation}
\xi\hat{*}\eta=\sum\limits_{k=0}^{n-1}\sum\limits_{i,j=0}^{n-1}
\hat{A}^k_{ij}\xi^i\eta^je_k
\end{equation}
where
\begin{equation}
\begin{array}{lcr}
\hat{A}_{ij}^k=\bar{A}_{ij}^k& \mbox{for}& 1\leq i,j,k\leq n-1 \\
\hat{A}_{ij}^0=\bar{g}_{ij}& \mbox{for}& 1\leq i,j\leq n-1 \\
\hat{A}_{0i}^j=0& \mbox{for}& 0\leq i,j\leq n-1 
\end{array}
\end{equation}
Next, if we factor $({\cal A}^n_{\hat{W}})^*$ over the ideal ${\bf C}e_0$ we
get what we formerly called coextension:
\begin{equation}
\begin{array}{c}
\xi\bar{*}\eta=\sum\limits_{k=1}^{n-1}(\sum\limits_{i,s,j=1}^{n-1}\xi^i\eta^j
\bar{A}_{ij}^k)e_k\\[8pt]
\xi,\eta\in ({\cal A}^{n-1}_{\bar{W}})^*
\end{array}
\end{equation}
If we want to come back, then naturally
\begin{Proposition}
The matrix $\left(\bar{g}_{ij}\right)_{1\leq i,j\leq n-1}$ defines an
extension of the algebra $({\cal A}^{n-1}_{\bar{W}})^*$:
\begin{equation}
\begin{array}{c}
e_i\hat{*}e_j=e_i\bar{*}e_j+\bar{g}_{ij}e_0=\sum\limits_{k=1}^{n-1}A_{ij}^k e_k+
\bar{g}_{ij}e_0\\[4pt]
e_i\hat{*}e_0=0\\[4pt]
e_0\hat{*}e_0=0\\[4pt]
(1\leq i,j\leq n-1)
\end{array}
\end{equation}
\end{Proposition}
The structure of $({\cal A}^{n+1}_W)^*$ is obtained simply adding the
semisimple part.

With the help of the new coextension structure  it is easy to
the equations for the Casimir components (\ref{eq:CC1}),
(\ref{eq:CC2}) together with $C_{0j}=0$ for $0\leq j\leq n-1$. They all can be
written in an elegant form which we shall introduce below. To this end 
we note that $\Psi$ satisfies the same sort of identity we had for
$\bar{\Psi}$:
\begin{equation}
T^*_{\Psi^{-1}(\xi)}\eta=\xi*\eta, \quad \xi,\eta\in ({\cal A}^n)^*
\end{equation}
Then we can simply write:
$$
C(e_i,e_j)=C(\Psi(\Psi^{-1}(e_i)),e_j)=C(T^*_{\Psi^{-1}(e_i)}e_n,ej)=
C(e_n,T^*_{\Psi^{-1}(e_i)}(e_j))
$$
and get the following elegant result:
\begin{equation}\label{eq:CC3}
C_{ij}=C(e_n,e_i*e_j),\quad 0\leq i,j\leq n-1
\end{equation}
or equivalently
\begin{equation}\label{eq:CC4}
C_{ij}=\sum\limits_{k=0}^{n-1}\hat{A}_{ij}^k C_{nk},\quad 0\leq i,j\leq n-1
\end{equation}
Note also that $e_i*e_j=e_i\hat{*}e_j$ for $1\leq i,j\leq n-1$, so if we
wish we can put $e_i\hat{*}e_j$ in the formula (\ref{eq:CC3}). As
$e_0*e_i=0$ the above formula not only gives more concise expression
than (\ref{eq:CC1}), (\ref{eq:CC2}) but also includes the equations $C_{i0}=0$.

\section{Conclusion}
We have considered the Casimir invariants, related to universal
extensions using the equivalence between these structures and
commutative algebra structures and have  shown that they have clearer
meaning for the corresponding commutative structures. In the
nondegenerate case we succeeded to give a clear algebraic meaning of the
so-called coextension structures.  We hope that this will stimulate
further studies of the Poisson-Lie structures of different physical
theories and their Casimir invariants.

\end{document}